\documentclass[12pt,a4paper]{article}
\usepackage{amsmath}
\usepackage{amsthm}
\usepackage{theorem}
\usepackage{makeidx}
\usepackage{amstext}
\usepackage[cp1251]{inputenc}
\usepackage{epsfig}
\begin{document}
\title{On numerical integration of coupled Korteweg-de
Vries System}
\author{A.A. Halim $^{1}$, S.P. Kshevetskii $^{2}$, S.B. Leble$^{3}$ \\
        $^{1,3}$ Technical University of Gdansk, \\
        ul. G. Narutowicza 11/12, 80-952 Gdansk, Poland. \\
        $^{2}$ Kaliningrad state University, Kaliningrad, Russia\\
         $^{3}$leble@mif.pg.gda.pl}
\maketitle
\begin{abstract}
We introduce a numerical method for general coupled Korteweg-de
Vries systems. The scheme is valid for solving Cauchy problems for
arbitrary number of equations with arbitrary constant
coefficients. The numerical scheme takes its legality by proving
its stability and convergence which gives the conditions and the
appropriate choice of the grid sizes. The method is applied to
Hirota-Satsuma (HS) system and compared with its known explicit
solution investigating the influence of initial conditions and
grid sizes on accuracy. We also illustrate the method to show the
effects of constants with a transition to non-integrable cases.
\end{abstract}

\section{Introduction}
Coupled Korteweg-de Vries system equations form a class of
important nonlinear evolution systems. Its importance comes
(physically) from the wide application field it covers and
(mathematically) from including both (weak) nonlinearity and third
order derivatives (weak dispersion). It describes the interactions
of long waves with different dispersion relations. Namely it is
connected with most types of long waves with weak dispersion
($\omega(k)\rightarrow0,(k)\rightarrow0$), e.g.  internal,
acoustic and planetary waves in geophysical hydrodynamics.\\

It was introduced by A.Maxworthy, L.Redekopp and P.Weldman
\cite{Max} in studying the nonlinear atmosphere Rossby waves.
R.Hirota  and J.Satsuma \cite{Hir} give single- and two-soliton
solutions to some version of the system. R.Dodd and A.Fordy
\cite{Dod} found an L-A pair for Hirota-Satsuma equations.
S.B.Leble  derived the cKdV system for different hydrodynamical
systems with explicit expressions for the nonlinear and dispersion
constants \cite{Leb}. He also developed the approach to the cKdV
integration. S.B.Leble, S.P.Kshevetskii \cite{Leb,Ksh} used the
system in investigation of nonlinear internal gravity waves. A.
Perelemova \cite{Per} used it in description of interaction of
acoustic waves with opposite directions of propagation
in liquids with bubbles.\\

Others deal with integrability of the system (Dodd R. and Fordy A.
\cite{Dod}) from a Lax pair point of view, S.B.Leble \cite{Leb2}
in Walquist-Estabrook theory. Foursov MV \cite{Four} described a
new method for constructing integrable system of differential
equations that reduced to cKdV equations. Oevel W. \cite{Oev}
considers integrable system of cKdV and found an infinite
hierarchy of commuting symmetries and conservation laws in
involution. Zharkov A Yu. \cite{Zhar} obtained a new class of
integrable KdV-like systems. Metin Gurses, Atalay Karasu
\cite{Gur} found infinitely many coupled system of KdV type
equations which are integrable. They also give recursion
operators. In studying the Painleve test classification of the
system, Ayse (Kalkanli) Karasu \cite{Ays} found new KdV systems
that are completely integrable in the sense of WTC paper. He was
looking for the integrable subclass of KdV systems given by
Svinolupov \cite{Svi}. The later has introduced a class of
integrable multicomponent KdV equations associated with Jordan
algebras. John Weiss \cite{Weis2} derived the associated
"modified" equations for HS system and from these the Lax pair is
also derived. B.A.Kupershmidt \cite{Kup} showed that a dispersive
system describing a vector multiplet interacting with the KdV
field is a member of a bi-Hamiltonian integrable hierarchy.\\

The significant achievement in numerical solution of the single
KdV equation starts from the famous paper of Norman Zabusky and
Martin Kruskal \cite{Z}. It develops the idea of soliton solutions
set for the integrable equations and enlighten the problems of
effective integration scheme elaborating. The paper launched the
numerous investigations and inventions in this field. Perhaps, the
last publication that develop applications of recent theoretical
achievements in numerical integration schemes is based on the
notion of isospectral deformations \cite{Cel}. Recently a
multisymplectic twelve points scheme was produced \cite{Ping}.
This scheme is equivalent to the multisymplectic Preissmann scheme
and is applied to solitary waves over long time interval.

Shaohong Zhu \cite{Shao} had produced a difference scheme for the
periodic initial-boundary problem of the Coupled KdV (Hirota-
Satsuma case) system. He use the inner product of the discrete
function to obtain a scheme keeping two conserved quantities. His
scheme is a nonlinear algebraic system for which a catch-ran
iterative method is designed to solve it.

The coupled KdV system representing most possible physical
application ( related to the weak nonlinear dispersion ) to be
considered in this work takes the following general form

\begin{equation}
\left( \theta _{n}\right) _{t}+c_{n}\left( \theta _{n}\right)
_{x}+\sum\limits_{k,m}g_{mkn}\,\theta _{k}\left( \theta
_{m}\right) _{x}+d_{n}\left( \theta _{n}\right)
_{xxx}=0\ , \,\,\,\  n,m,k=1,2,3,...,N , \\
\end{equation}
where $\theta _{n}$ (x, t) is the  amplitude of the wave mode as a
function of space x and time t respectively. The constants $c_{n}$
are the linear velocities and \thinspace $ g_{mkn},\,e_{n}$ are
the nonlinear and dispersion coefficients.\\

In the present work we introduce a numerical tool for solving
coupled KdV system which is a development of the two step three
time levels as Lax-Wendroff scheme \cite{Ksh,Tan}. Proving the
theorem about stability and convergence of the scheme gives the
opportunity to use it for different applications like Cauchy
problems for arbitrary number of equations and a wide class of
initial conditions $\theta _{n}$ (x, 0). We consider in our
problem an infinite domain while the initial condition goes
quickly enough to zero following the relation
\begin{equation*}
\int(1+|x|)|\theta(x,0)|dx=0<\infty
\end{equation*}
keeping in mind the choice of smooth and integrable function. As
an important corollary of the theorem one obtains conditions that
have to be taken in account in choosing grid sizes. This numerical
method is checked by applying it to HS system for which a good
number of explicit solutions exist \cite{Leb3}. We examine also
the effects of equations coefficients and conditions of the
problem on the solution.

In the section 2 we introduce the difference scheme for arbitrary
number of coupled KdV equations. We investigate stability and
prove the convergence giving the condition have to be taken in
account in choosing the grid sizes and how they are related. In
section 3 we analyze HS system with two-parameters one-soliton
explicit solution. The numerical method is applied to HS system
and compared with the explicit solution. We analyzed the effects
of the two parameters and initial condition on the form of the
resulting solitons as will as on accuracy and show the results by
figures. We also produced (numerically) a multi-soliton solution
for HS system and used the conservation law to estimate the
expected number of solitons which agreed that we already obtained.
Proving stability and convergence besides testing the results for
HS system allows us in section 4 to use scheme for general cKdV
system. Hence we illustrate by plots the results of applying the
scheme to slightly nonintegrable cKdV systems and other for a
system with non-smooth initial conditions.

\section{The numerical method}
\subsection{The difference scheme}
For the cKdV system (1) we introduce a numerical (
finite-difference ) method of solution. A scheme which is two
steps three time levels similar to the Lax-Wendroff one
\cite{Ksh,Tan}. The usual Lax-Wendroff is modified such that the
order of the first derivative becomes of order $O(\triangle
x^{4})$. The approximation of the nonlinear terms is changed in
such a manner that the integral of $\theta ^{2}$ be a conserved
one. The approach gives a solution that can be considered as some
generalized solution, in the sense of Shwartz distribution theory,
where the dispersion constants vanishes. This scheme is suitable
to a nonlinear equations and is valid for n equations with
arbitrary coefficients. The scheme can be simply derived beginning
from Taylor series expansion as
\begin{equation}
\left( \theta _{n} \right) _{i}^{j+1} =\left( \theta _{n}
\right)_{i}^{j} + \Delta t\,\left(\left( \theta _{n}
\right)_{t}\right)^{j}_{i} +O\lbrack (\Delta t)^{2} ].
\end{equation}
where i and j are used to locate a point in the discrete domain
and $\Delta t$ is the time step while the subscript t means
time derivative. Substituting for  $%
 \left( \theta _{n} \right)_{t}  $ in (2) using the
system (1) to obtain
\begin{equation}
\left( \theta _{n} \right) _{i}^{j+1} =\left( \theta _{n} \right)
_{i}^{j} - \Delta t \left( c_{n} \left( \theta _{n} \right) _{x} +
\sum\limits_{k,m}g_{mkn} \theta _{k} \,\left( \theta _{m} \right)
_{x} + d_{n} \left( \theta _{n} \right) _{xxx} \right)^{j}_{i}
+O\,\,\left[ \left( \Delta \,t\right) ^{2} \right].
\end{equation}
The difference scheme is elaborated applying Lax idea for a half
time step and leap frog method to the remaining half time step. In
both steps $\left( \theta _{n} \right) _{x} and\,\left( \theta
_{n} \right) _{xxx} $  are replaced by forth order $O\left( \Delta
x^{4}\right)$ and second order accurate $O\left( \Delta
x^{2}\right)$ central difference expressions. Hence (3) gives the
following difference scheme
\begin{multline}
\left(\left( \theta _{n} \right) _{i}^{j+\frac{1}{2} } -\left(
\theta _{n} \right) _{i}^{j} \right)/{\frac{\tau }{2} }
+c_{n}\left({\left( \theta _{n} \right) _{i+1}^{j} -\left( \theta
_{n} \right) _{i-1}^{j} }\right)/{2h}+ \sum\limits_{k,m}g_{mkn}
\left( \theta _{k} \right) _{i}^{j} \left({\left( \theta _{m}
\right) _{i+1}^{j} -\left( \theta _{m} \right) _{i-1}^{j}
}\right)/{2h} \tag{4.1}
\\
+e_{n} \left({\left( \theta _{n} \right) _{i+2}^{j} -2\left(
\theta _{n} \right) _{i+1}^{j} +2\left( \theta _{n} \right)
_{i-1}^{j} -\left( \theta _{n} \right) _{i-2}^{j} }\right)/{ 2
h^{3} } =0 , e_{n}=\left(d_{n} -{c_{n}h^{2}}/{6}\right) \notag,
\end{multline}
where n, m, k are the modes numbers; i and j are discrete space
and time respectively. The time step $\Delta t$ is replaced for
simplicity by t while h denotes spatial step. The equation (4.1)
is accompanied with discrete equation for the intermediate layer
as:\\
\begin{multline}
\left({\left( \theta _{n} \right) _{i}^{j+1} -\left( \theta _{n}
\right)_{i}^{j} }\right)/{\tau }
+c_{n} \left({\left( \theta _{n} \right) _{i+1}^{j+\frac{1}{%
2} } -\left( \theta _{n} \right) _{i-1}^{j+\frac{1}{2} }
}\right)/{2h} +\sum\limits_{k,m}g_{mkn} \left( \theta _{k} \right)
_{i}^{j+\frac{1}{2} } \left({\left( \theta _{m} \right)
_{i+1}^{j+\frac{1}{2} } -\left( \theta _{m} \right)
_{i-1}^{j+\frac{1}{2} } }\right)/{2h} \tag{4.2}
\\
+e_{n} \left({\left( \theta _{n} \right) _{i+2}^{j+\frac{1}{2} }
-2\left( \theta _{n} \right) _{i+1}^{j+\frac{1}{2} } +2\left(
\theta _{n} \right) _{i-1}^{j+\frac{1}{2} } -\left( \theta _{n}
\right) _{i-2}^{j+\frac{1}{2} } }\right)/{ 2 h^{3} } =0 \notag
\end{multline}

\subsection{Stability and convergence analysis}

For simplicity of the analysis we start by considering one
equation of the system and give  the details of stability and
convergence. Then we apply the idea to the general cKdV system
because it is rather close to that for one KdV equation but
  more bulky.

\subsubsection{Stability analysis for KdV scheme}

Consider one KdV equation of the system (1)
\begin{equation}
\theta _{t}^{} +c_{} \theta _{x}^{} +g_{} \,\theta ^{} \theta
_{x}^{} +d_{} \theta _{xxx}^{} =0 \tag{5}
\end{equation}
Note again that the investigation we perform can be generalized
for the case of any finite number of modes.
Considering the numerical scheme applied for the equation (5)\\
\begin{equation}
(\theta _{i}^{j+1}-\theta _{i}^{j})/\tau +c(\theta
_{i+1}^{j}-\theta _{i-1}^{j})/2h+g\theta _{i}^{j}(\theta
_{i+1}^{j}-\theta _{i-1}^{j})/{2h}+e(\theta _{i+2}^{j}-2\theta
_{i+1}^{j}+2\theta _{i-1}^{j}-\theta _{i-2}^{j})/ 2 h^{3}=0
\tag{6}
\end{equation}

Let us select a suitable norm. For this multiply equation
equation (5)
 by $\theta$ and integrate to yield\\
\begin{equation*}
\frac{1}{2} \frac{d}{ dt} \int\nolimits_{-\infty }^{\infty }\theta
^{2} dx =0 \,\,\,\,\ or \,\,\,\,\,\int\nolimits_{-\infty }^{\infty
}\theta ^{2} dx = const,
\end{equation*}
hence, by definition of $L_{2} $  norm,$\left( \left\| \theta
\right\| _{2} \right)^{2}=\int\nolimits_{-\infty }^{\infty }\theta
^{2} dx   $ ,it may be written as \,\,\,\,\,\,\,\,\ $\left(
\left\| \theta \right\| _{2} \right) ^{2} =\,\,$ const., i.e. the
norm $\left\| \theta \right\| _{2} $ is
conserved and the equation can be treated in the $L_{2} $  norm.\\

Now we will prove stability with respect to small perturbations (
because we consider nonlinear equations ) of initial conditions.
Strictly speaking it is the boundness of the discrete solution in
terms of small perturbation of the initial data. So let us
consider the
differential\\
\begin{equation}
d\theta _{i}^{j+1} =\sum\limits_{r} ({\partial \,\theta _{i}^{j+1}
}/{\partial \,\theta _{r}^{j} }) d\theta _{r}^{j} \,\ ,\,\,\,
\,\,\, r=...,i-1,i,i+1,...\tag{7}
\end{equation}
for equation (5) denoting
 $\,\,T_{i,r}^{j+1} =\left\{ {\partial
\,\theta _{i}^{j+1} }/{\partial \,\theta _{r}^{j} } \right\}
\,\,\,\,\, ,\,\,\,\,\,\,d\theta _{r}^{j} =\left\{
\begin{array}{l}
d\theta _{i-2}^{j} \\
d\theta _{i-1}^{j} \\
d\theta _{i}^{j} \\
d\theta _{i+1}^{j} \\
d\theta _{i+2}^{j}
\end{array}
\right\} $ and use $\left\| d\theta ^{j}\right\| =\left(
\sum\limits_{r}\left( d\theta _{r}^{j}\right) ^{2}\,\, h  \right) ^{\frac{1}{2}}$%
.\\
Rewrite also the relation (7) in the matrix form\\
\begin{equation*}
d{\theta}_{i}^{j+1} =T^{j+1}_{i,r} d{\theta}^{j}_{r} =T^{j+1}T^{j}
d{\theta}^{j-1} =\prod\limits_{r}T^{r} d{\theta}^{0},
\end{equation*}
 where $d{\theta}^{0}$ is a small perturbation of
initial data and the subscripts are omitted for simplicity.\\

Stability requires the boundedness of the product \,\,\,\
$\prod\limits_{r}T^{r} $ \,\,\ in a sense that the norm $\left\|
\prod\limits_{r}T^{r}\right\| $ is bounded by some constant, i.e.
$\left\| \prod\limits_{r}T^{r}\right\| \leq C$. Here $C$ is a
constant, and the matrix norm is a spectral norm. For this the
sufficient condition is
 \ $\left\| T^{r}\right\| <e^{a\tau }$ where $a$ is a constant, that
is independent of $\tau $.
 The case $\left\|T^{r}\right\|<e^{a(t,h)\ast \tau }$ is a sufficient
condition of stability also, but only if $\left| a(\tau
,h)\right| \leq const<\infty $. If  $\left| a(\tau ,h)\right|
<const$, including $\tau $, $h\rightarrow 0$, for some dependence
$\tau =f(h)$, then we can say about conditional stability. Namely
this kind of stability will be climbed below.\\

To calculate $T^{r}$, rewrite the scheme (6) in the form
\,\,\,\,\,\,\ $\theta _{i}^{j+1}=\theta _{i}^{j+1}(\theta
_{i+2}^{j},\theta _{i+1}^{j},\theta _{i}^{j},\theta
_{i-1}^{j},\theta _{i-2}^{j}).$
 \,\,\,\,\\
So
\begin{multline}
\left(T ^{j+1}\right)_{ir} =\delta _{i,r} -({c \tau }/{2h})
\left[ \delta _{i+1,r} -\delta _{i-1,r} \right] -({g\tau
}/{2h})\left[ \theta _{i}^{j} \left(\delta _{i+1,r} -\delta
_{i-1,r} \right) +\delta _{i,r} \left( \theta _{i+1}^{j} -\theta
_{i-1}^{j}\right)\right] \tag{8}
\\
-(e\tau /{2h^{3} }) \left[ \delta _{i+2,r} -2\delta _{i+1,r}
+2\delta _{i-1,r} -\delta _{i-2,r} \right]. \notag
\end{multline}
 Rewriting (8) in terms of the
identity (E), symmetric (S) and anti-symmetric (A) matrices\\
$T^{j+1}=C+S^{j+1}+A^{j+1}$
\begin{equation*}
\left\{S^{j+1}\right\}_{i,r}=-\frac{g\tau }{4h}\left(\left(\theta
_{i}^{j}-\theta _{i+1}^{j}\,\right) \delta _{i+1,r}-\left( \theta
_{i}^{j}-\theta _{i-1}^{j}\right) \delta _{i-1,r}+2\delta
_{i,r}\left[ \theta _{i+1}^{j}-\theta _{i-1}^{j}\right] \right)
\end{equation*}
\begin{eqnarray*}
\left\{ A^{j+1}\right\} _{i,r} &=&-\frac{c\tau }{2h}\left[ \delta
_{i+1,r}-\delta _{i-1,r}\right] -\frac{g\tau }{4h}\,\left( \left(
\theta _{i}^{j}+\theta _{i+1}^{j}\right) \,\delta _{i+1,r}-\left(
\theta_{i}^{j}+\theta _{i-1}^{j}\right) \,\delta _{i-1,r}\right) \\
&&-\frac{e\tau }{2h^{3}}\left[ \delta _{i+2,r}-2\delta
_{i+1,r}+2\delta _{i-1,r}-\delta _{i-2,r}\right]
\end{eqnarray*}
\begin{equation*}
\left\| S^{j+1}\right\| \leq \left| g\right| \tau
\,\max_{i}(\left| \theta _{x,i}^{j}\right|,\left| \theta\grave{~}
_{x,i}^{j}\right|) \qquad ,\theta _{x,i}^{j}=\left[ \theta
_{i+1}^{j}-\theta _{i}^{j}\right]/(h)\,\ and \,\ \theta\grave{~}
_{x,i}^{j}=\left[ \theta_{i+1}^{j}-\theta _{i-1}^{j}\right]/(2h)
\end{equation*}\\
one arrive at
\begin{equation*}
\left\| A^{j+1}\right\| \leq \frac{\left| g\right| \tau
}{h}\,\max_{i}\left| \theta _{i}^{j}\right| +\frac{\left| c\right|
\tau }{h}\,\,\,+ \frac{3\left| e\right|\tau }{h^{3}}.
\end{equation*}

\begin{eqnarray*}
\left\| T^{j+1}\right\|^{2}&=&\left\| \left( T^{j+1}\right)^{\ast
}T^{j+1}\right\| =\left\| \left( E-A^{j+1}+S^{j+1}\right) \left(
E+A^{j+1}+S^{j+1}\right) \right\| \\
&\leq& 1+2\left\| S^{j+1}\right\| +\left( \left\| A^{j+1}\right\|
+\left\|S^{j+1}\right\| \right) ^{2} \\
&\leq& 1+2\left| g\right| \tau \,\max_{i}\left| \theta
_{x,i}^{j}\right|+ \tau ^{2}\left( 2\left| g\right|
\,\max_{i}\left| \theta _{x,i}^{j}\right| +\frac{\left| g\right|
}{h}\,\max_{i}\left| \theta _{i}^{j}\right| +\frac{\left| c\right|
}{h}\,\,\,+\frac{3\left| e\right| }{h^{3}}\right) ^{2} \\
&\leq& \,\,\,\,e^{a\tau }\,\,\,\,\,where\,\,a=2\left| g\right|
\,\max \left| \theta _{x,i}^{j}\right| +\tau \left( 2\left|
g\right| \,\max_{i}\left| \theta _{x,i}^{j}\right| +\frac{\left|
g\right| }{h} \,\max_{i}\left| \theta _{i}^{j}\right|
+\frac{\left| c\right| }{h}\,\,\,+ \frac{3\left| e\right|
}{h^{3}}\right) ^{2}\\
\end{eqnarray*}
which is a necessary condition of stability. The scheme is stable
if $a\leq constant\leq \infty $ in spite of $\tau ,$ $h\rightarrow
0.$ This is a conditional stability of the scheme. It means that
it is required for stability that $\tau \rightarrow 0$ more
faster than $h\rightarrow 0$,\,\ or

\begin{equation}
\tau \leq (constant). \,\,\ h^{6} , \,\,\,\ constant<\infty
\tag{9}
\end{equation}

Therefore, for small enough $\tau $ we can simplify the
expression for a

\begin{equation*}
a=2\left| g\right| \,\max \left| \theta _{x,i}^{j}\right| +\tau
\left( \frac{3e}{h^{3}}\right) ^{2}
\end{equation*}
In practical calculations the time step $\tau $ should be chosen
so that it would satisfy\,\ $\tau \left( \frac{3e}{h^{3}}\right)
^{2}\ast t_{0}=O(1)$, where $t_{0}$ is the time of simulation
($0\leq t\leq t_{0}$). In future, when we shall be suggesting some
better numerical scheme, we will essentially use our observation
that stability depends only on the dispersion terms. And now we
will try to accomplish our short investigation of the scheme (6)
by strong proving of the numerical scheme convergence.

\subsubsection{ The proof of the KdV scheme convergence}

Now we prove that a solution of equation (6) converges to a
solution of (5), if the exact solution is a
continuously-differentiable one. Let us denote by $\theta (x,t)$ a
solution of the equation (5). We substitute $\theta _{i}^{j}=$
$\theta (x_{i},t_{j})+v_{i}^{j}$ into (6), \ $v_{i}^{j}$ is a
error between the difference solution $\theta _{i}^{j}$ and the
exact solution $\theta (x_{i},t_{j})$. We obtain the equation for
$v_{i}^{j}$
\begin{multline}
\left(v_{i}^{j+1}-v_{i}^{j}\right)/\tau+c
\left(v_{i+1}^{j}-v_{i-1}^{j}\right)/2h +g
\theta(x_{i},t_{j})\left(v_{i+1}^{j}-v_{i-1}^{j}\right)/2h +g
v_{i}^{j}\left( \theta (x_{i+1},t_{j})-\theta
(x_{i-1},t_{j})\right)/2h\\
+g v_{i}^{j}\left(v_{i+1}^{j}-v_{i-1}^{j}\right)/2h
+e\left(v_{i+2}^{j}-
2 v_{i+1}^{j}+ 2 v_{i-1}^{j}- v_{i-2}^{j}\right)/ 2 h^{3}= \\
-\left(\left(\theta(x_{i},t_{j+1})-\theta(x_{i},t_{j})\right)/\tau+c\left(
\theta (x_{i+1},t_{j})-\theta (x_{i-1},t_{j})\right)/2h+ g
\theta(x_{i},t_{j})\left(\theta (x_{i+1},t_{j})- \theta
(x_{i-1},t_{j})\right)/2h\right.\\
\left. +e \left(\theta (x_{i+2},t_{j})-2\theta
(x_{i+1},t_{j})+2\theta (x_{i-1},t_{j})-\theta (x_{i-2},t_{j})
\right)/  2 h^{3} \right) \notag
\end{multline}
Let us take into account that
\begin{multline*}
v_{i}^{j}-\tau \left( c(v_{i+1}^{j}-v_{i-1}^{j})/2h+g\theta
(x_{i},t_{j})(v_{i+1}^{j}-v_{i-1}^{j})/2h+gv_{i}^{j}(\theta
(x_{i+1},t_{j})-\theta (x_{i-1},t_{j}))/2h\right.  \\
\left. +e(v_{i+2}^{j}-2v_{i+1}^{j}+2v_{i-1}^{j}-v_{i-2}^{j})/ 2
h^{3}\right) =\sum_{k}\left( T^{j+1}\right) _{ik}v_{k}^{j}
\end{multline*}
Using the operator $T^{j+1}$ introduced above, this equation may
be rewritten in the form
\begin{multline}
\left( v_{i}^{j+1}-\sum_{k}\left(
T^{j+1}\right)_{ik}v_{k}^{j}\right)/\tau+g
v_{i}^{j}\left(v_{i+1}^{j}-v_{i-1}^{j}\right)/2h =\\
-\left(\left(\theta(x_{i},t_{j+1})-\theta(x_{i},t_{j})\right)/\tau+c
\left( \theta (x_{i+1},t_{j})-\theta (x_{i-1},t_{j})\right)/2h+ g
\theta(x_{i},t_{j})\left(\theta (x_{i+1},t_{j})- \theta
(x_{i-1},t_{j})\right)/2h\right.\\
\left. +e \left(\theta (x_{i+2},t_{j})-2\theta
(x_{i+1},t_{j})+2\theta (x_{i-1},t_{j})-\theta (x_{i-2},t_{j})
\right)/ 2 h ^{3} \right) \notag
\end{multline}
The right part of this relation is a quantity of order $O(\tau
+h^{2}).$ So, we can write
\begin{equation}
\left( v_{i}^{j+1}-\sum_{k}\left( T^{j+1}\right)
_{ik}v_{k}^{j}\right) /\tau +g
v_{i}^{j}\left(v_{i+1}^{j}-v_{i-1}^{j}\right)/2h=O(\tau +h^{2}),
\notag
\end{equation}
or
\begin{equation}
v_{i}^{j+1} =\sum_{k}\left( T^{j+1}\right) _{ik}v_{k}^{j}-\tau
f_{i}^{j}, \,\,\,\,\,\,\ f_{i}^{j}
=gv_{i}^{j}\frac{v_{i+1}^{j}-v_{i-1}^{j}}{2h}-O(\tau +h^{2})
\tag{10}.
\end{equation}
We finally arrive at the inequality that compare the norms:
\begin{equation*} \left\| f^{j}\right\| \leq \frac{\left|
g\right| }{h^{\frac{3}{2}}}\left\| v^{j}\right\| ^{2}+O(\tau
+h^{2})
\end{equation*}
To explain how this estimate was obtained, follow the expressions
\begin{multline}
\left\| f^{j}\right\|  =\left( \sum\limits_{i}\left(
f_{i}^{j}\right) ^{2}h\right) ^{\frac{1}{2}} \leq \left| g\right|
\left( \sum\limits_{i}\left(
v_{i}^{j}\frac{v_{i+1}^{j}-v_{i-1}^{j}}{2h}\right) ^{2}h\right)
^{\frac{1}{2}}+O(\tau +h^{2}) \notag \\
\leq \left| g\right| \left( \sum\limits_{i}\left( v_{i}^{j}\right)
^{2}h\ast \sum\limits_{i}\left( v_{i}^{j}\right) ^{2}h\right)
^{\frac{1}{2}} \frac{1}{h^{\frac{3}{2}}}+O(\tau +h^{2}). \tag{11}
\end{multline}
Using the Schwartz inequality $\left\| AB\right\| \leq \left\|
A\right\| \left\| B\right\| $, the formulas (10) could be
transformed as
\begin{multline}
\left\| v^{j+1}\right\| \leq \left\| T^{j+1}\right\| \left\|
v^{j}\right\| +\tau \left\| f^{j}\right\| \leq \left\|
T^{j+1}\right\| \left\| T^{j}\right\| \left\| v^{j-1}\right\|
+\tau (\left\| T^{j+1}\right\| \left\|
f^{j-1}\right\| +\left\| f^{j}\right\| )\leq   \tag{12} \\
\left\| T^{j+1}\right\| \left\| T^{j}\right\| \left\|
T^{j-1}\right\| \left\| v^{j-2}\right\| +\tau (\left\|
T^{j+1}\right\| \left\| T^{j}\right\| \left\| f^{j-2}\right\|
+\left\| T^{j+1}\right\| \left\| f^{j-1}\right\|
+\left\| f^{j}\right\| \leq   \notag \\
e^{a\tau j}\left\| v^{0}\right\| +\tau (e^{a\tau (j-1)}\left\|
f^{0}\right\| +e^{a\tau (j-2)}\left\| f^{1}\right\| +...\left\|
f^{j}\right\| )\leq
\notag \\
e^{a\tau j}\left\| v^{0}\right\| + M\max_{k\leq j}(\left\|
f^{k}\right\|
)\leq   \notag \\
e^{a\tau j}\left\| v^{0}\right\| +M\left( \frac{\left| g\right| }{h^{\frac{3%
}{2}}}\left\| v^{j+1}\right\| ^{2}+O(\tau +h^{2})\right) ,\qquad
M=\tau \frac{e^{a\tau j}-1}{e^{a\tau }-1}.  \notag
\end{multline}
To derive (12), we  have \bigskip used the iteration of the first
of formula (we substituted the formula into itself, but for index
less than 1) and using (11). Then we have utilized the formula for
a sum of geometric series. Farther the inequality obtained in (12)
has a solution
\begin{equation*}
\left\| v^{j+1}\right\| \leq \left(1-\sqrt{1-4M\frac{\left|
g\right| }{h^{
\frac{3}{2}}}\left( e^{a\tau }\left\| v^{0}\right\| +MO(\tau +h^{2})\right) }%
\right)/\left(2M\frac{\left| g\right| }{h^{\frac{3}{2}}}\right)
\end{equation*}
If we take into account (9), and use$\left\|v^{0}\right\|=0,$we
obtain
\begin{multline}
\left\| v^{j+1}\right\|  \leq \left(1-\sqrt{1-\frac{4M\left|
g\right| }{h^{ \frac{3}{2}}}MO(\tau
+h^{2})}\right)/\left(\frac{2M\left|
g\right| }{h^{\frac{3}{2}}}\right)\\
\approx \left(1-\left(1-\frac{2M\left| g\right| }{h^{
\frac{3}{2}}}MO(\tau +h^{2})\right)\right)/\left(\frac{2M\left|
g\right| }{h^{\frac{3}{2}}}\right)=
 M\ast O(\tau +h^{2}) \notag
\end{multline}
The constant $M$ is bounded, in spite of $j\rightarrow \infty $, because of $%
j\tau <\infty $. Therefore, the convergence is proved.

\subsubsection{The coupled KdV scheme }
The numerical scheme for the system of the cKdV (4.1)-(4.2) is
also conditionally stable and convergent one. The  proof for this
scheme is close to that given before, but a bit bulky. We deal
with a vector:
\begin{equation*}
U =\left\{
\begin{array}{l}
\theta _{1}\,\,\ \theta _{2} \,\,\ ... \,\,\,\ \theta _{N}
\\\end{array}
\right\}^{t}
\end{equation*}
as a dependent variable instead of the simple variable $\theta $
in the case of one KdV. For this vector case the norm used has the
form
\begin{equation*}
\left\| U \right\|=
\left(\sum\limits_{l=1}^{N}\sum\limits_{i}\left| \theta
_{l,i}\right| ^{2}\,\, h \right) ^{\frac{1}{2}}.
\end{equation*}
The conditions connecting time and space steps for this scheme
look also similar but with different constants
\begin{equation*}
\max_{n}(\left| e_{n}\right| )\ast \frac{81\ast \tau ^{3}}{4\ast
h^{12}}\ast t_{0}=O(1).
\end{equation*}
\section{Checking the numerical method}
The numerical method is tested by applying it to Hirota-Satsuma
system. Namely the two parameters one-soliton explicit solution is
used \cite{Leb3}.
\subsection{Analytic solution ( explicit formula ) of Hirota-Satsuma
system} Darboux transformation (DT) that account a deep reduction
for this specific HS case of cKdV is used \cite{Leb3} to obtain
explicit solutions to HS system. The Lax representation of the HS
equations is based on the matrix 2x2 spectral problem of the
second order. For this problem the deep reduction scheme
\cite{Leb3} is applied (with the help of the conserved bilinear -
forms) and supports the constrains on the potential while the
iterated DT are performed. The iterated DT in determinant form and
the covariance of the bilinear forms with respect to DT under
restrictions gives N soliton solution of HS system. The system of
HS we use to check the scheme has the form:
\begin{equation}
\left( \theta _{1} \right) _{t} \,-0.25\,\left( \theta _{1}
\right) _{3x} \,-1.5\,\left( \theta _{1} \right) _{x} \,\left(
\theta _{1} \right) {}+3\left( \theta _{2} \right) _{x} \,\left(
\theta _{2} \right) =0 \tag{13.a}
\end{equation}
\begin{equation}
\left( \theta _{2} \right) _{t} \,+0.5\,\left( \theta _{2}
\right)_{3x}+1.5\,\left( \theta _{2} \right) _{x} \,\left( \theta
_{1} \right) \,=0. \tag{13.b}
\end{equation}
This system has two-parameters one soliton solution:
\begin{equation}
\begin{array}{cc}
\theta_{1} \,={-2m^{2} (-1+d^{2}
+2d\,Sin(\lambda_{1})*Sinh(\lambda_{2}))}/{\left(
d\,Cos(\lambda_{1})+Cosh(\lambda_{2})\right) ^{2} }, \tag{14} \\
\theta _{2} ={ \left( 2+2 d^{2} \right) ^{.5} m^{2}
}/{\left( d\,Cos(\lambda_{1})+Cosh(\lambda_{2})\right) },\\
\lambda_{1}=.5m^{3}t+mx \,\,\,\,\ and \,\,\,\,\
\lambda_{2}=.5m^{3}t-mx.
\end{array}
\end{equation}
with real constants m, d.  For small $\left| d\right| $ this
solution is a smooth function but for $\left| d\right| \,>1$ poles
appear .The following figures show some choice of m and d to show
the effect of these two parameters on the solution. Figure 1 show
that, for constant d, the amplitude is proportional to m while the
wave width is inversely proportional to it. Figure 2 show that,
for the given m, d affects on soliton shape, namely the first
mode, while the amplitude of the second is inversely proportional
to d.
\begin{center} \epsfig{file=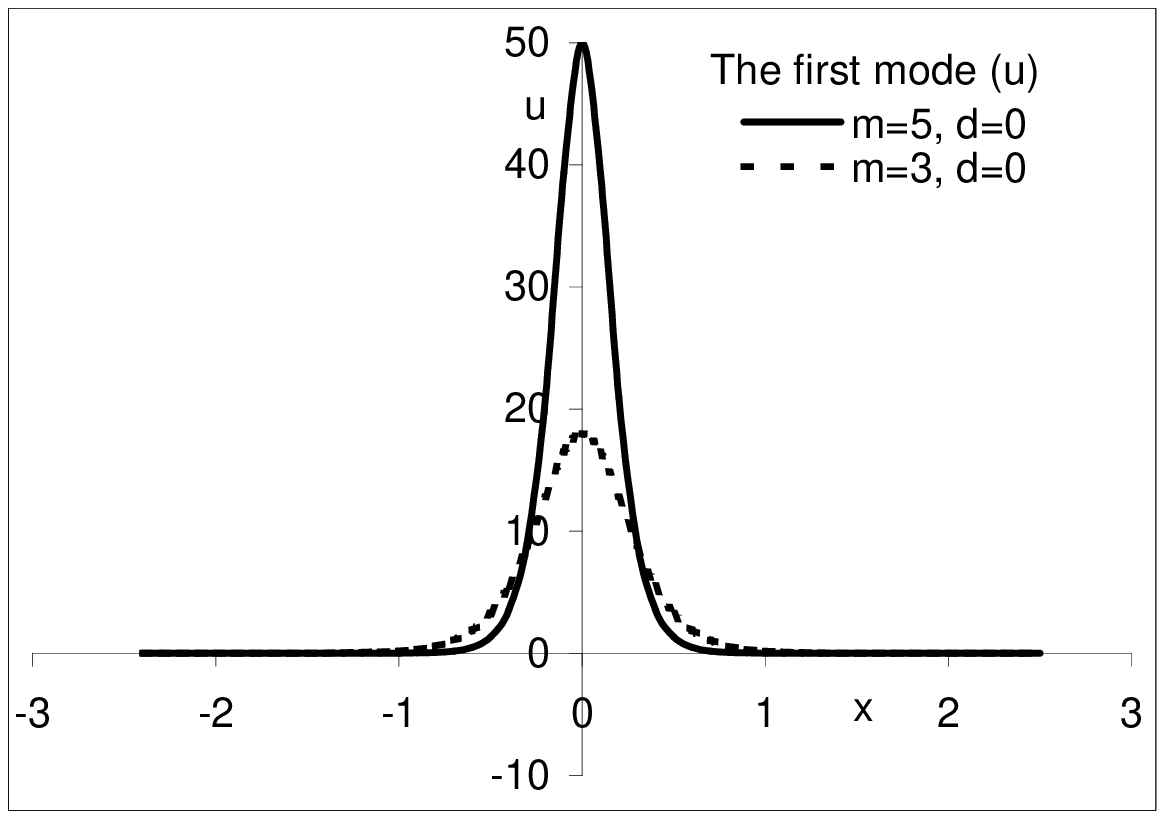, height=3.7cm,
width=8.6cm ,clip=,angle=0} \epsfig{file=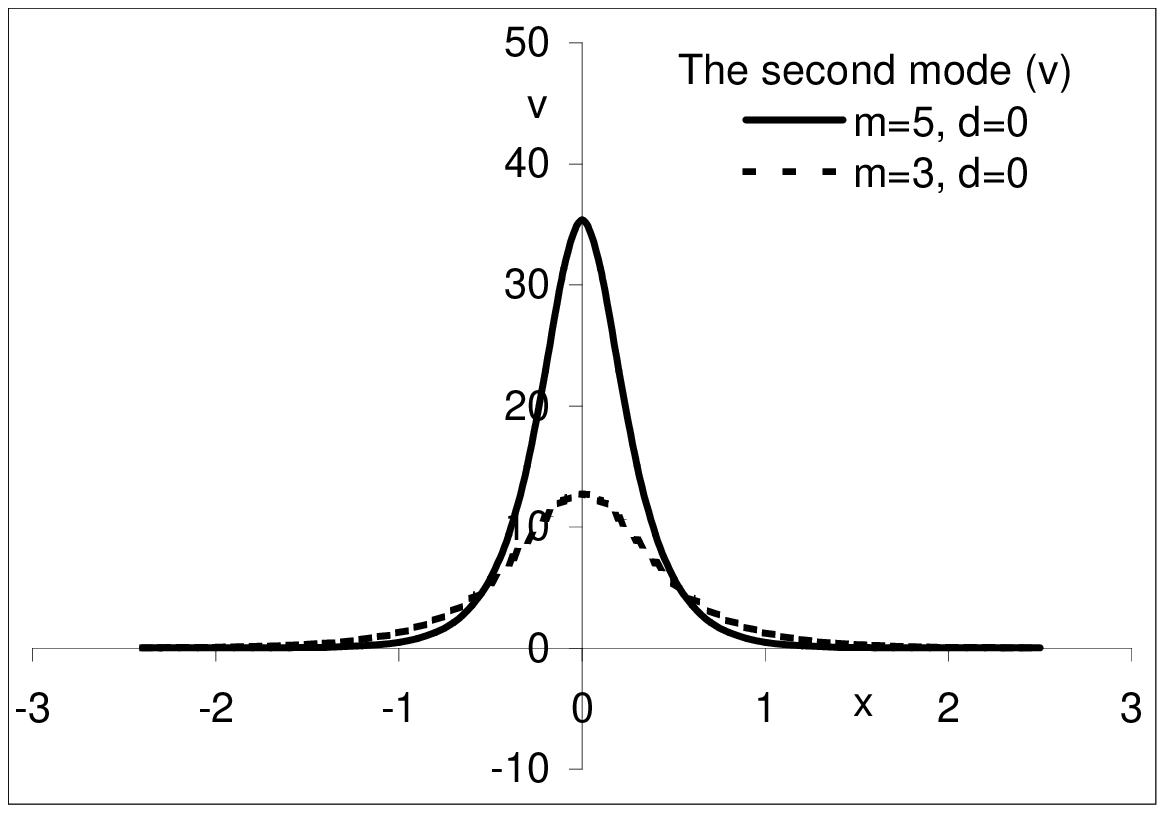,
height=3.7cm, width=8.6cm ,clip=,angle=0}\\
\ Fig.1. For a
constant d, the amplitude is proportional to m while the wave
width is inversely proportional to m.
\end{center}
\begin{center}
\epsfig{file=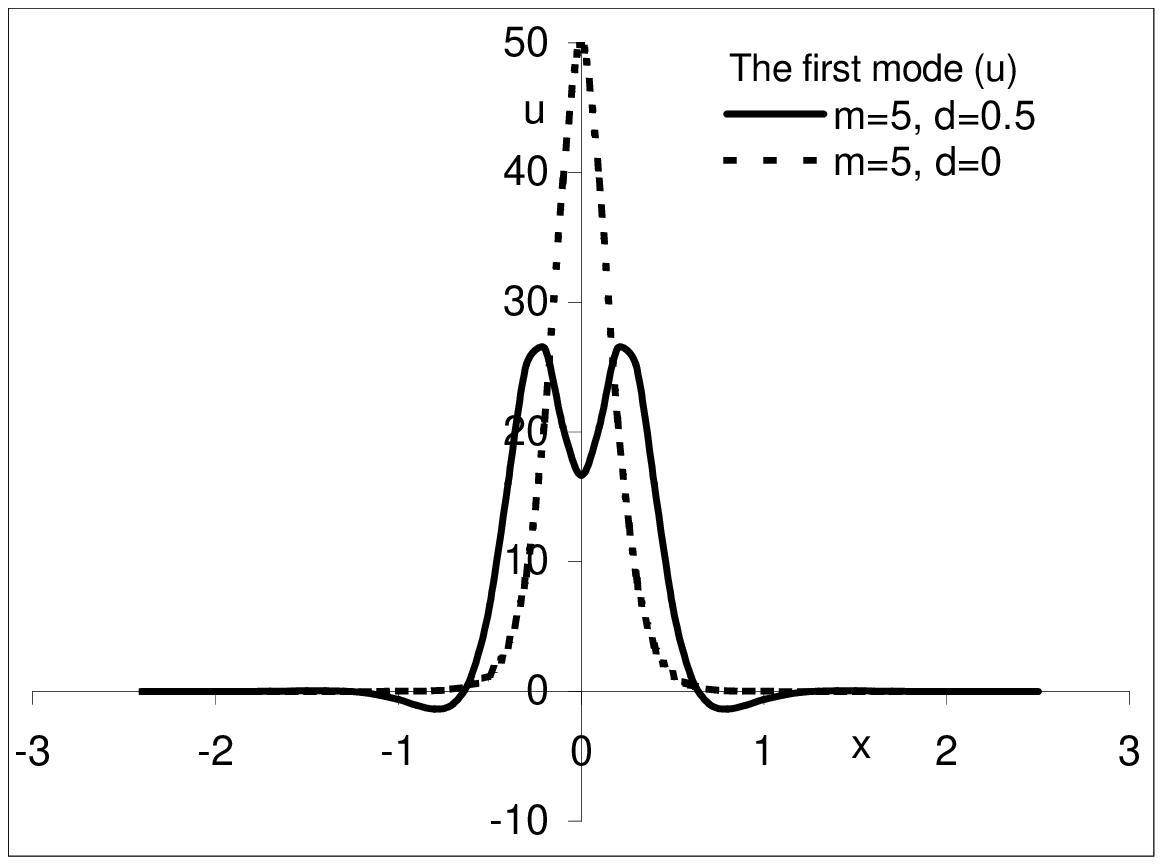, height=3.7cm, width=8.6cm
,clip=,angle=0} \epsfig{file=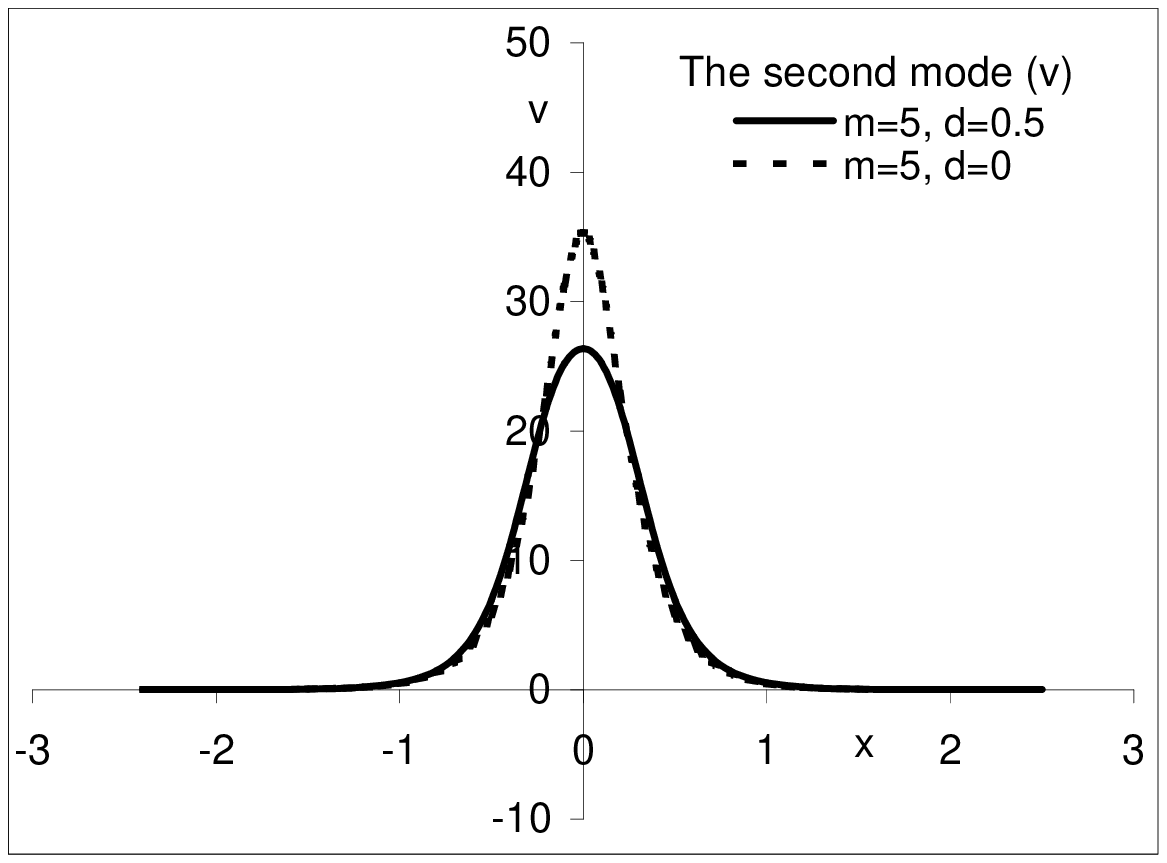, height=3.7cm,
width=8.6cm ,clip=,angle=0}\\
 \ Fig. 2. For the same m, d affects
on soliton shape, namely the first mode, while the amplitude \\ of
the second mode is inversely proportional to d.
\end{center}
\subsection{Calculations by numerical scheme and comparison results}
HS system (13) is solved numerically using the scheme (4) with
initial condition from formula (14)(t = 0) and the results are
compared with the explicit solution. It is found that, keeping the
restriction on the choice of $\tau \,\ and\,\,\ h $ and relation
between them, the initial wave modes amplitude affects on the
accuracy of the results . Also the error decreases as the mesh is
refined. Namely smaller amplitude ( of order one ) gives better
results as shown below in figure (3). It gives the percentage
error calculated as follow
\begin{equation}
\% Error=  \frac{|Explicit\,\ solution - Numerical \,\
solution|}{Initial\,\ amplitude } \notag
\end{equation}
We relate the error to the initial amplitude to show the physical
significance of the error.
\begin{center}
\epsfig{file=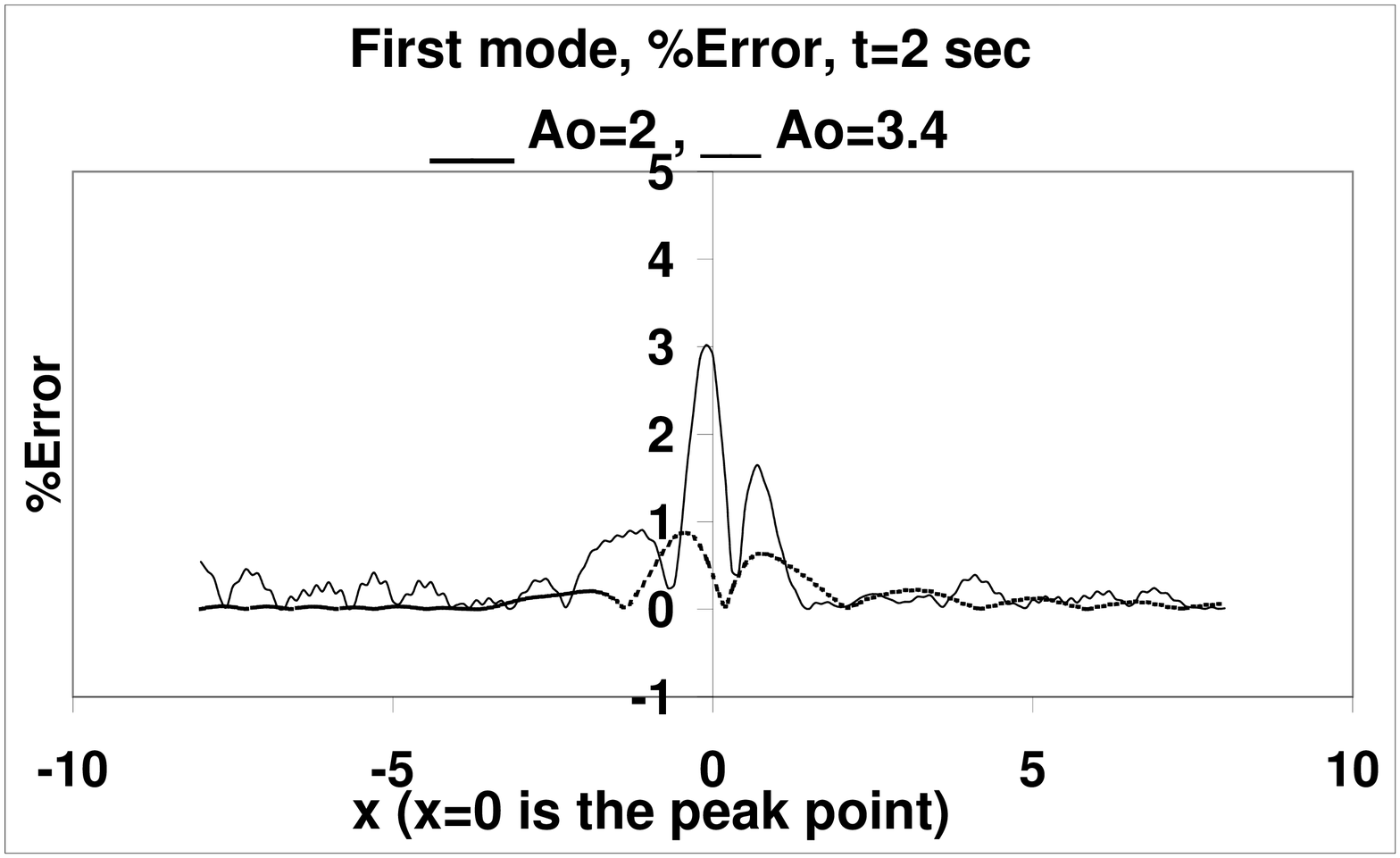, height=4cm, width=8.6cm ,clip=,angle=0}
\epsfig{file=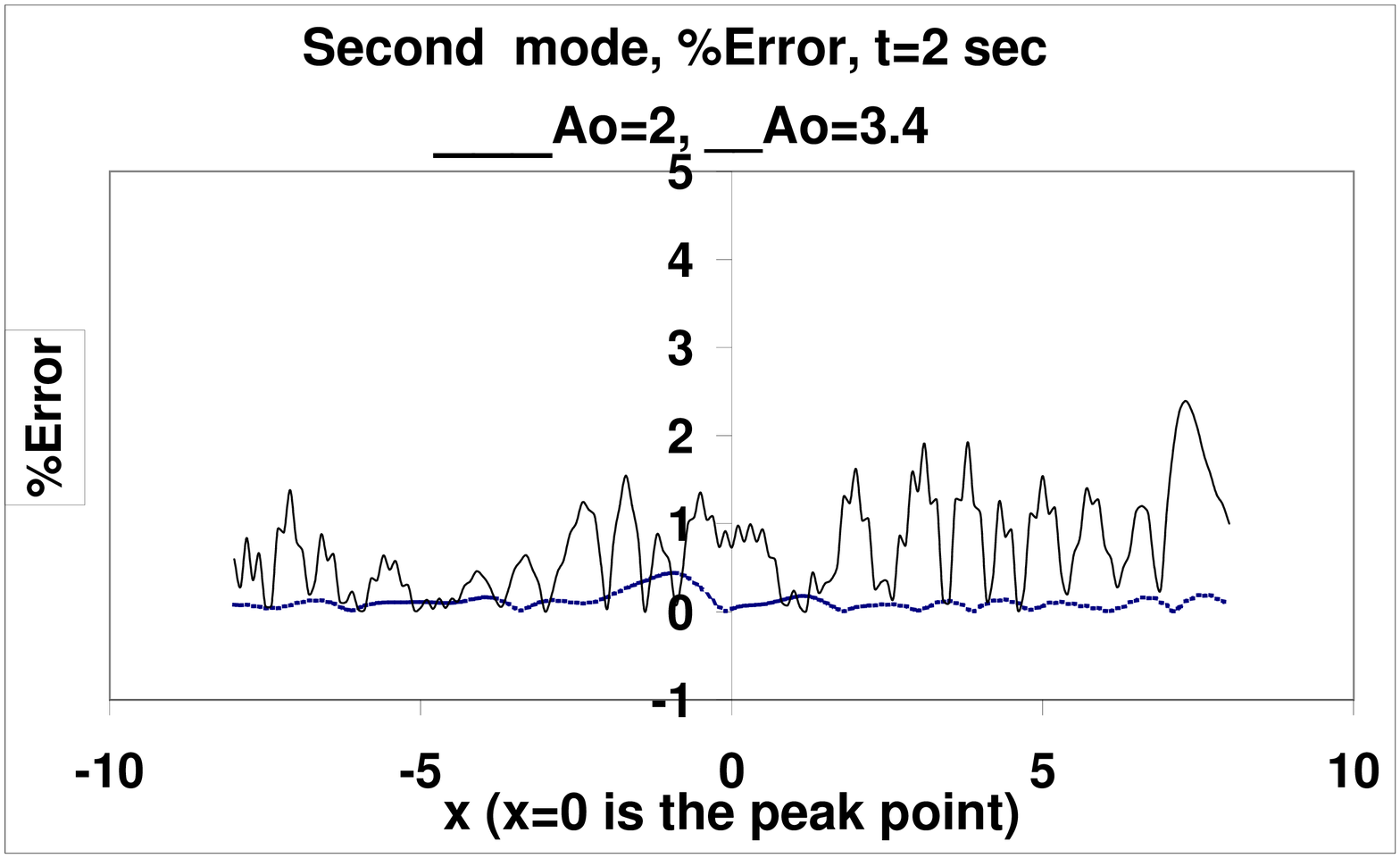, height=4cm, width=8.6cm ,clip=,angle=0}\\
{\hspace*{10 mm}  \ Fig.3 }\% error is proportional to the
amplitude(A).
\end{center}
The plots show that the error is proportional to the amplitude
(A), where as shown in figure the maximum relative error in the
case (A=2) is 1 \% while in the case (A=3.4) is 3 \%. The reason
may be due to the higher velocity in the larger one, hence more
interactions impact. It also shows that the error increases near
the peak points. The reason of these  oscillations in plots
appearance is that the numerical and analytical plots intersect
over the space domain.
\section{Applying the scheme to different applications}
Stability analysis and checking performed to the scheme in the
general cKdV equations  give  the ability of using this scheme to
solve other problems for which analytical solutions have not been
found. We first consider the multi-soliton solution decay of the
localized initial condition for the single KdV equation of HS
system (figure 4.a) then for the complete HS system (figure 4.b).
In both we use the initial condition from formula (14) but with 10
times the width and twice the amplitude.
\begin{center}
\epsfig{file=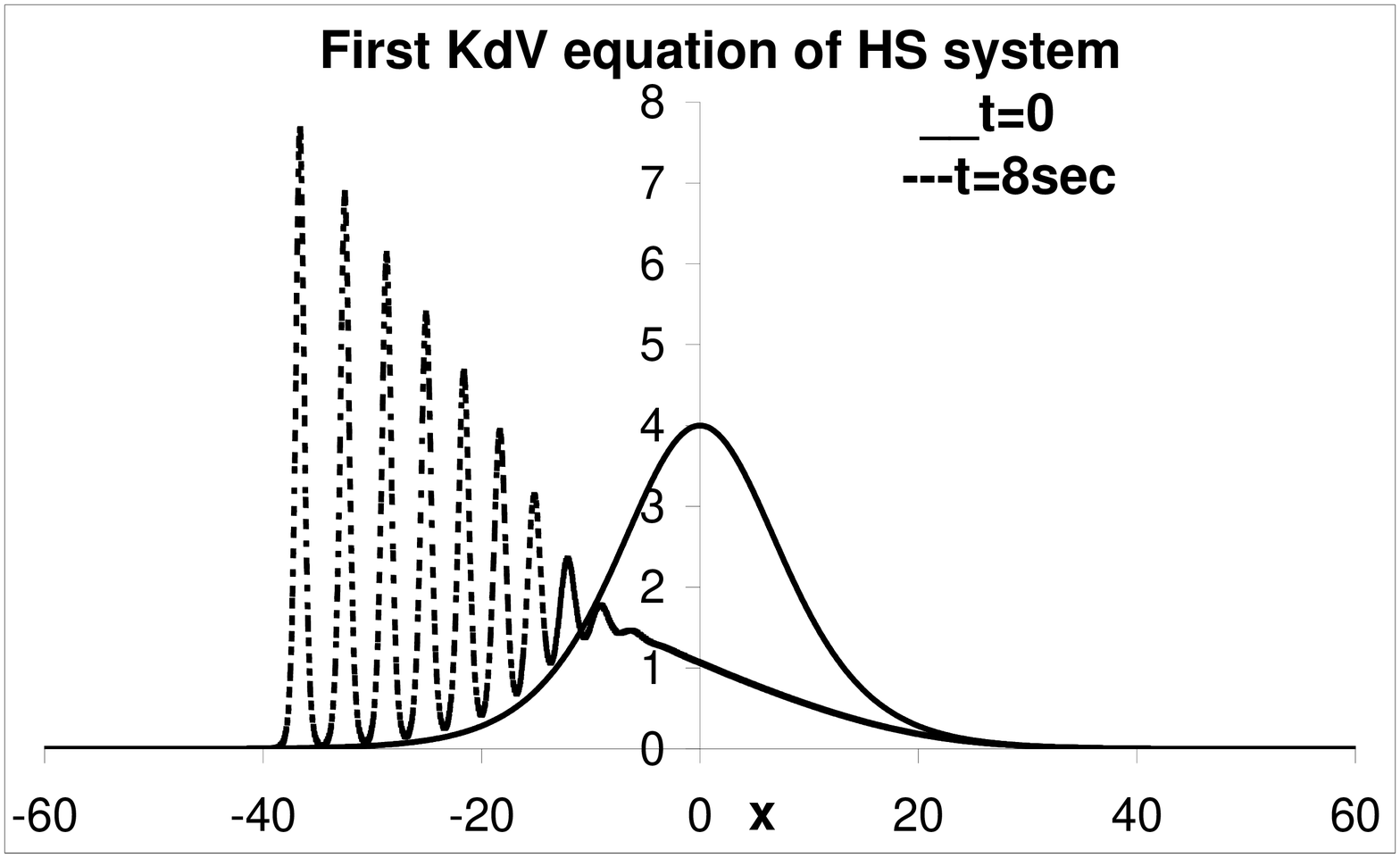, height=3.5cm, width=10cm
,clip=,angle=0}\\
Fig. 4.A  The multi-soliton decaying of the isolated (first) KdV
of HS system (13.a)
\end{center}
\begin{center}
\epsfig{file=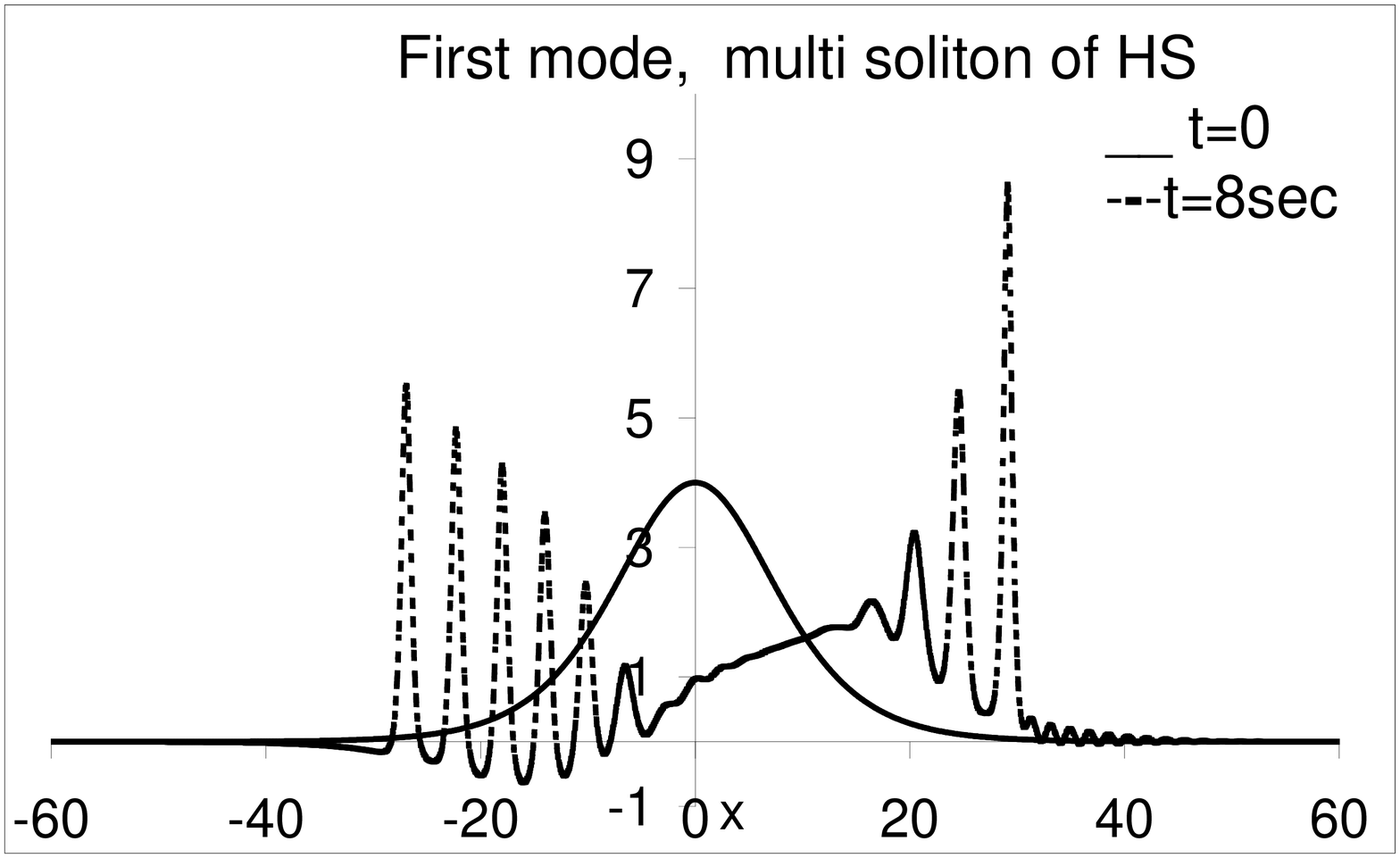, height=4cm, width=8.6cm ,clip=,angle=0}
\epsfig{file=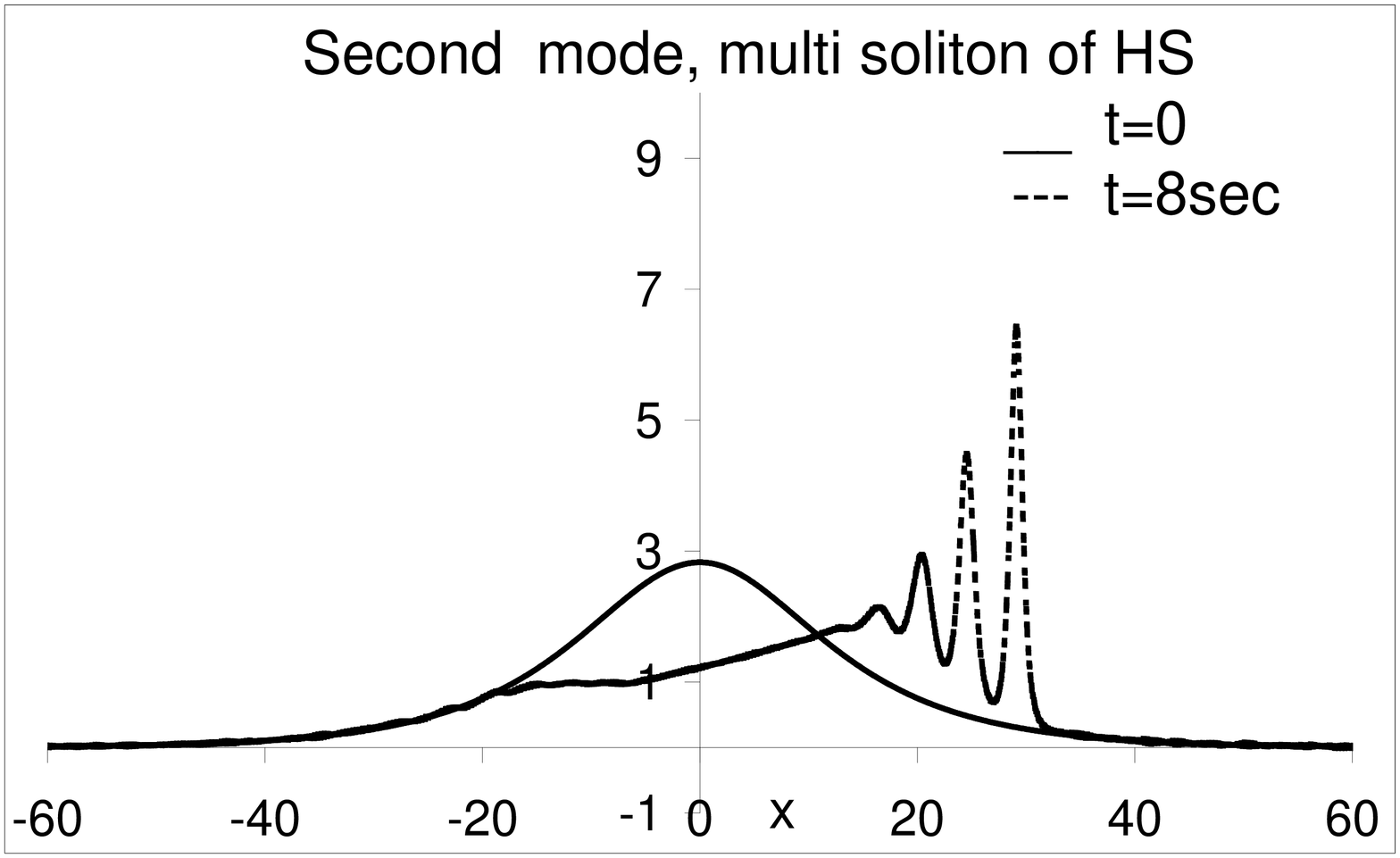, height=4cm, width=8.6cm
,clip=,angle=0}\\ Fig. 4.b The multi-soliton decaying of the
complete HS system (13).
\end{center}
The second mode affects (interacts) the first one which results in
the right direction soliton-like "tail" as shown in the figure
4.b. We estimate, using the conservation low (derived below), the
expected number of
solitons that already obtained in the numerical solution as\\
\begin{equation}
\begin{array}{cc}
\theta_{1} \times (13.1)-2\theta_{2}\times (13.2)\Rightarrow \notag \\
\frac{d}{dt}[.5 \theta_{1}^{2}-\theta_{2}^{2}]+\frac{d}{dx}[-.5
\theta_{1}^{3}-.25
\theta_{1}\theta_{1xx}+.125\theta_{1x}^{2}-\theta_{2}\theta_{2xx}
+.5\theta_{2x}^{2}]=0 \\
\theta_{1}\,\ and \,\ \theta_{2} \rightarrow 0 \,\ as \,\ x
\rightarrow \pm \infty\\
\int^{-\infty}_{\infty}{(.5
\theta_{1}^{2}-\theta_{2}^{2})dx}=const.
\end{array}
\end{equation}
Next we go to the solution of nonintegrable HS system. The
integrable HS system (13.a,b) may be shifted to "slightly"
nonintegrable one by small change of the dispersion constant of
the first equation to have the new nonintegrable HS system
\begin{equation}
\left( \theta _{1} \right) _{t} \,-0.2\,\left( \theta _{1} \right)
_{3x} \,-1.5\,\left( \theta _{1} \right) _{x} \,\left( \theta _{1}
\right) {}+3\left( \theta _{2} \right) _{x} \,\left( \theta _{2}
\right) =0 \tag{15.a}
\end{equation}
\begin{equation}
\left( \theta _{2} \right) _{t} \,+0.5\,\left( \theta _{2}
\right)_{3x}+1.5\,\left( \theta _{2} \right) _{x} \,\left( \theta
_{1} \right) \,=0. \tag{15.b}
\end{equation}
Using our scheme with initial condition from (14) we find that the
scheme works satisfactory (in the sence of convergence) even for
nonintegrable HS system as shown from Figure 5 below. The solution
looks like a soliton one for small time.
\begin{center}
\epsfig{file=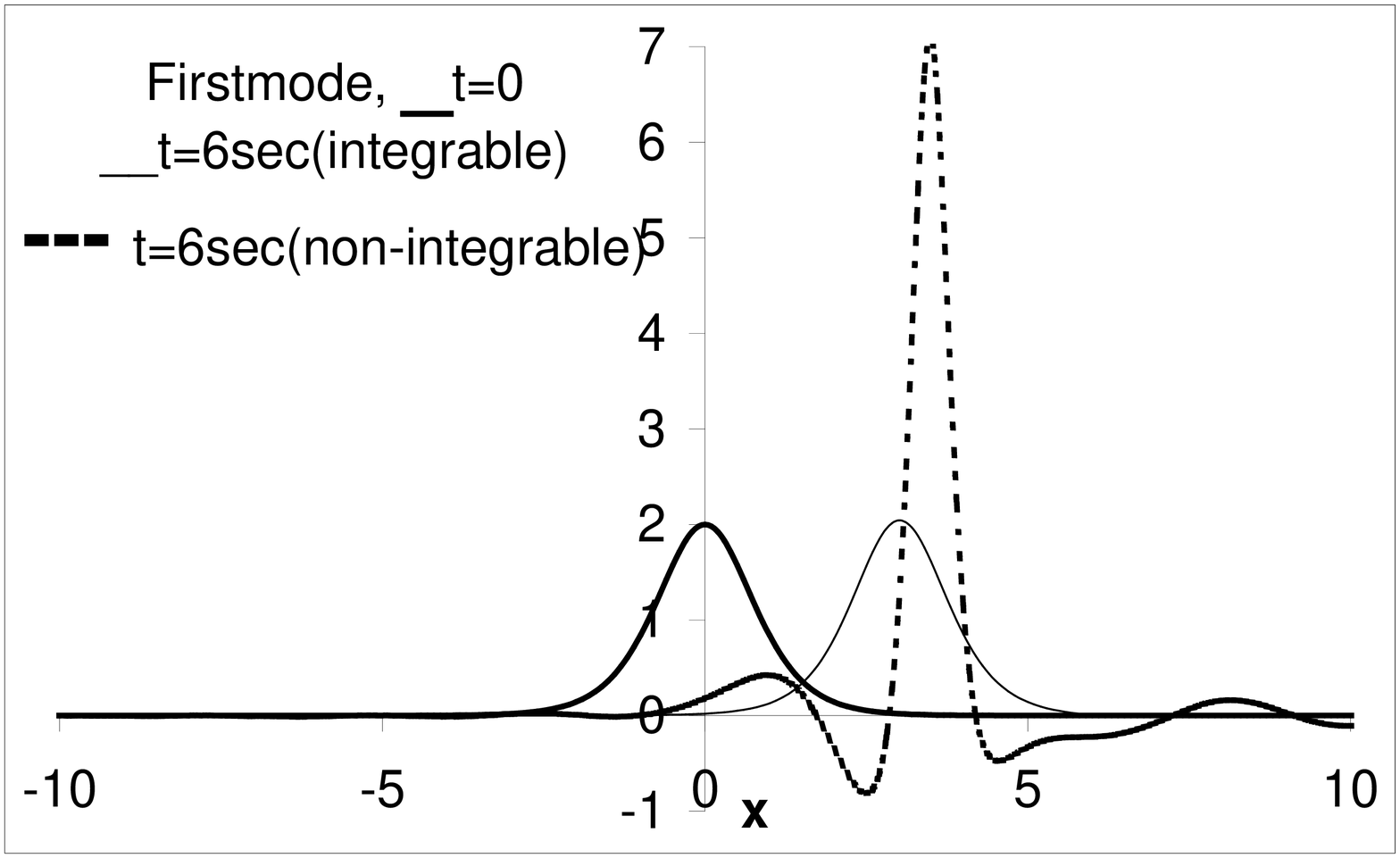, height=4cm,width=8.6cm,clip=,angle=0}
\epsfig{file=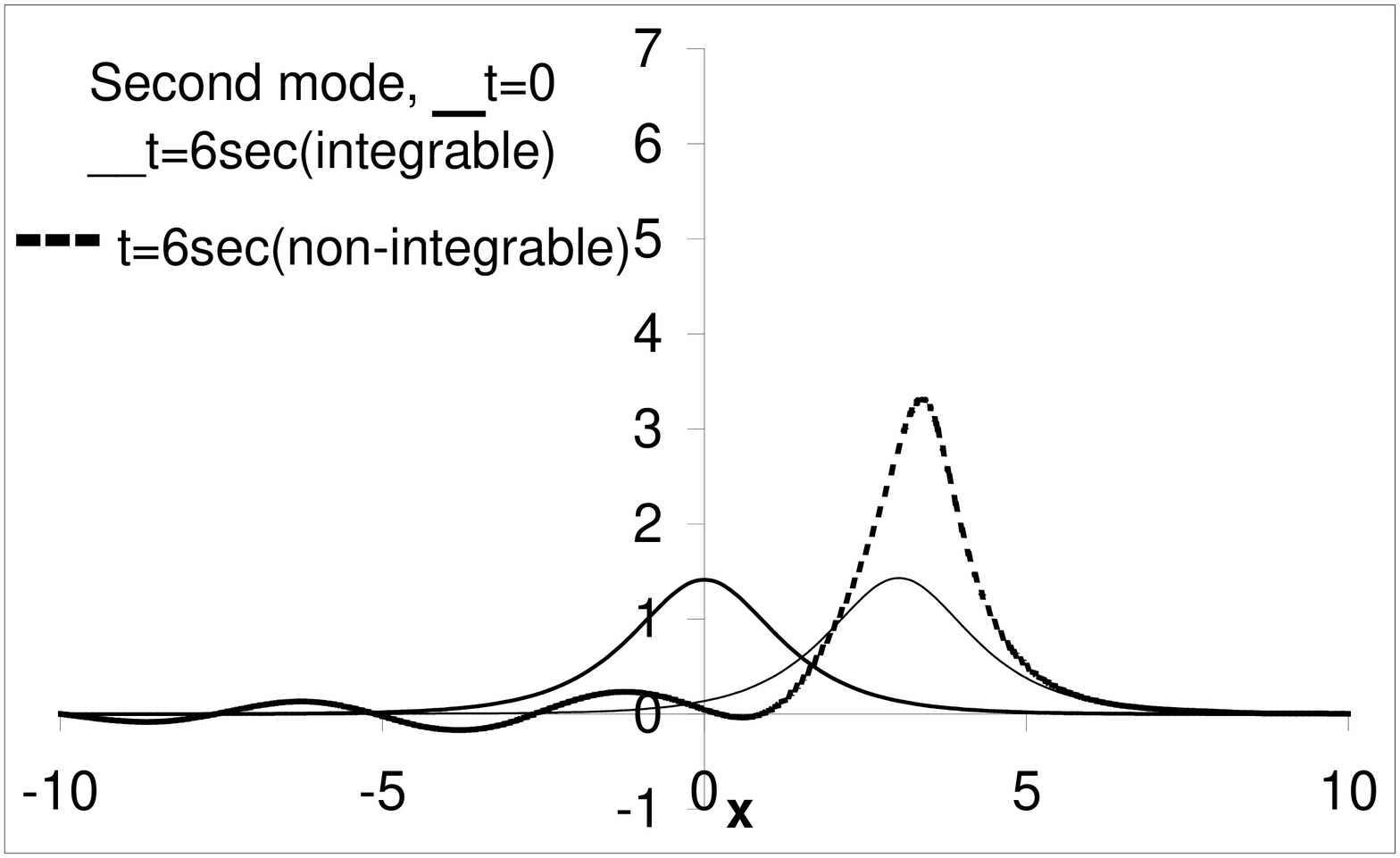, height=4cm,width=8.6cm ,clip=,angle=0}\\
\ Fig.(5) The numerical solution of integrable and slightly
nonintegrable HS system.
\end{center}
Also the solution using a non-smooth initial condition for HS
system (13) is shown in figure 6 below .
\begin{center}
\epsfig{file=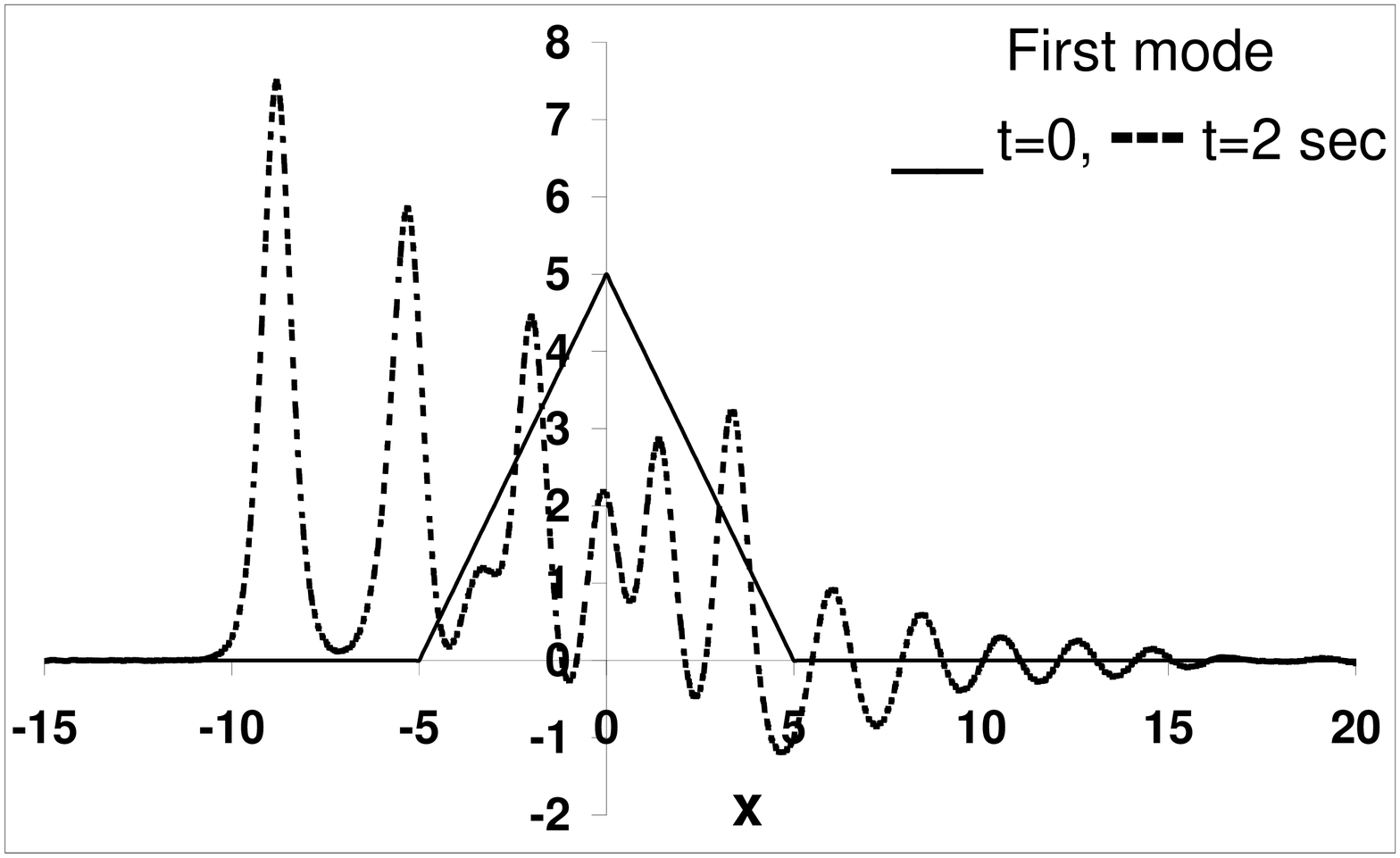, height=3.5cm,width=8.6cm,clip=,angle=0}
\epsfig{file=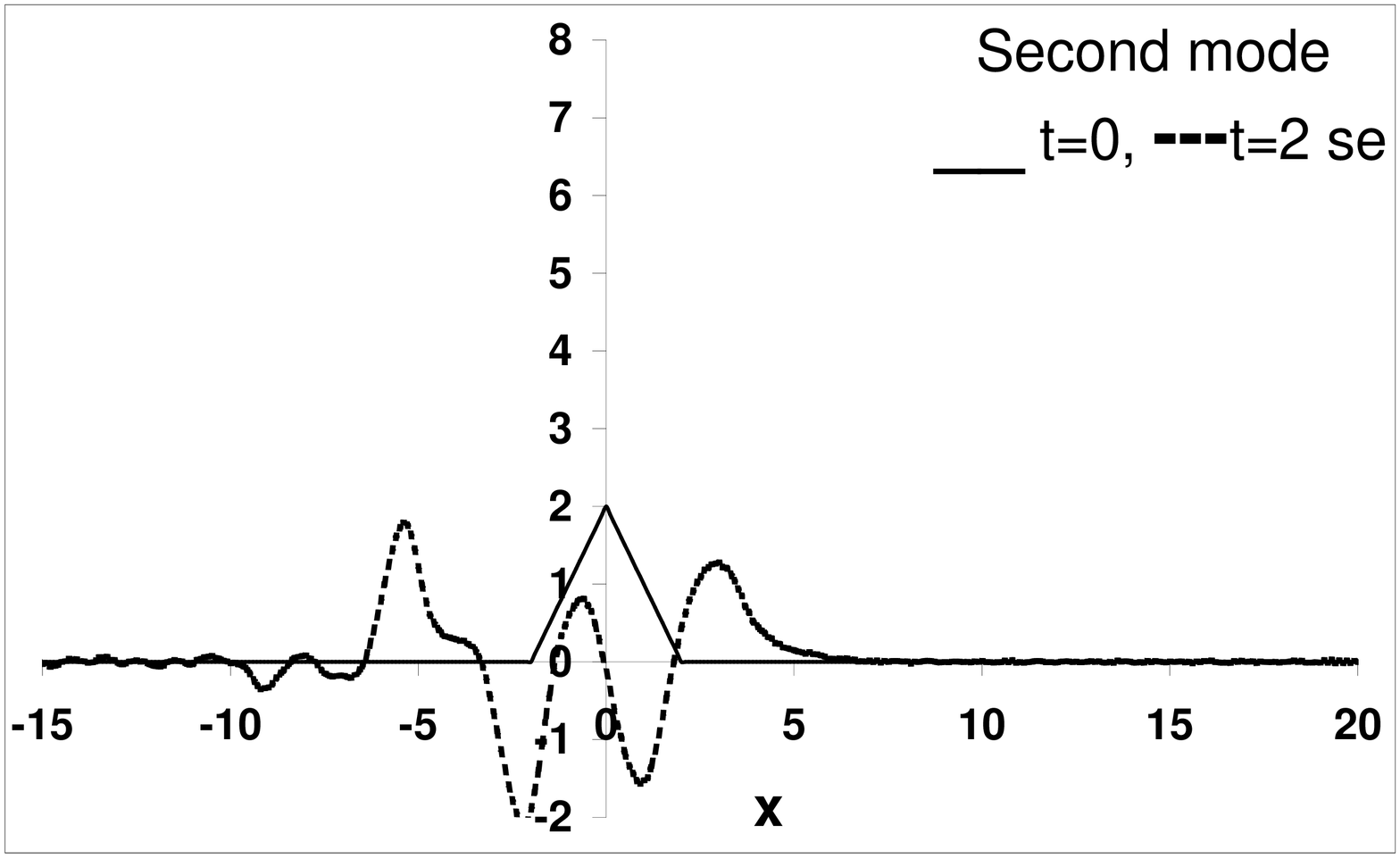, height=3.5cm,width=8.6cm ,clip=,angle=0}\\
\ Fig.6 The numerical solution of Non-smooth initial condition for
HS system (13).
\end{center}

\end{document}